\documentclass[a4paper,12pt]{article}
\usepackage{color}
\usepackage{amsmath}
\usepackage{graphicx}
\usepackage[authoryear]{natbib}

\usepackage{pdfcomment}

\usepackage{hyperref}

\usepackage{anziamjedraft}
\title{Macroscale boundary conditions for a non-linear heat exchanger}

\author{Chen Chen}
\address{School of Mathematical Sciences, University of Adelaide,
South Australia~5005, \textsc{Australia}.}
\mailto{chen.chen@adelaide.edu.au}

\author{A. J. Roberts}
\address{School of Mathematical Sciences, University of Adelaide,
South Australia~5005, \textsc{Australia}.}
\mailto{anthony.roberts@adelaide.edu.au}

\author{J. E. Bunder}
\address{School of Mathematical Sciences, University of Adelaide,
South Australia~5005, \textsc{Australia}.}
\mailto{Judith.Bunder@adelaide.edu.au}

\date{\today}

\usepackage{microtype,color}
\usepackage{graphicx} 

\usepackage{pgfplots}
\usepackage{tikz}
\usepackage{pgf}

\usetikzlibrary{arrows,external}

\newcommand{\dd}[2]{\frac{\partial #1}{\partial #2}}
\newcommand{\DD}[2]{\frac{\partial^2 #1}{\partial #2^2}}

\newcommand{\Ord}{\mathcal{O}}

\hypersetup{bookmarks,dvips,plainpages=false,
    colorlinks,linkcolor=red,citecolor=red,
    urlcolor=magenta,filecolor=magenta,
    breaklinks,bookmarksopen,pagebackref
}

\begin{document}

\maketitle

\begin{abstract}
Multiscale modelling methodologies build macroscale models of materials with complicated fine microscale structure.  We propose a methodology to derive boundary conditions for the macroscale model of a prototypical non-linear heat exchanger.  The derived macroscale boundary conditions improve the accuracy of macroscale model.  We verify the new boundary conditions by numerical methods.  The techniques developed here can be adapted to a wide range of multiscale reaction-diffusion-advection systems.
\end{abstract}

\tableofcontents

\section{Introduction\label{sec:introctac}}

Multiscale modelling techniques are a developing area of
research in engineering and physical sciences. These techniques are
needed when the system being modelled possesses very different
space-time scales and it is infeasible to simulate the whole domain
on a microscale mesh \citep{dolbow2004multiscale, Kevrekidis:2009fk, Bunder2012}.  Macroscale boundary conditions are rarely derived systematically in such works.
Instead macroscale boundary condition are often proposed heuristically  \citep{PS08, MR2777986, multiscalethesis}.
We developed a systematic method to derive boundary condition for
one-dimensional linear problems with fine structure
by cell mapping \citep{Chen2014}. Here we extend the method to a prototypical non-linear heat exchanger problem.

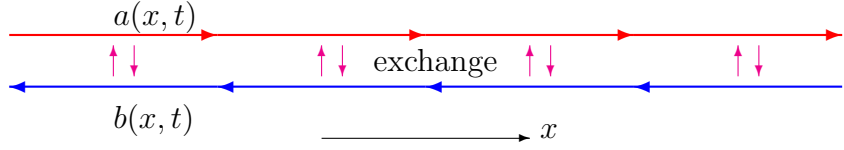
\begin{figure}
 \centering\setlength{\unitlength}{0.01\linewidth}
\begin{picture}(80,14)
\put(-10,0){
\put(20,1){\(b(x,t)\)}
\put(20,11){\(a(x,t)\)}
\put(45,6.5){exchange}
\multiput(20,0)(20,0)4{\color{magenta}
\put(0,6){\vector(0,+1){3}}
\put(2,9){\vector(0,-1){3}}
}
\put(40,0){\vector(1,0){20}\ $x$}
\thicklines
\multiput(10,10)(+20,0)4{\color{red}\vector(+1,0){20}}
\multiput(90,5)(-20,0)4{\color{blue}\vector(-1,0){20}}
}
\end{picture}
\caption{A schematic diagram of a heat exchanger. The red pipeline carries fluid to the right, and the blue pipeline carries fluid to the left. Heat exchanges between the pipes.\label{fig:he} }
\end{figure}

We mathematically model the counter
flow two-stream heat transfer shown in Figure~\ref{fig:he}.
Let \(x\)~measure nondimensional distance along the heat exchanger which is of length~\(L\), \(0\leq x\leq L\)\,, and let \(t\)~denote nondimensional time.
The field~$a(x,t)$ is the temperature of the fluid in one pipe and field~$b(x,t)$ is that in the other pipe. A quadratic reaction is included as an example nonlinearity to give the nondimensional microscale \text{pde}s
\begin{subequations}\label{eq:heatExchanger}%
\begin{eqnarray}
\dd{a}{t} & = & \tfrac{1}{2}\left(b-a\right)+\tfrac{1}{2}a^{2}-\dd{a}{x}+3\dd{^{2}a}{x^{2}}\,, \\
\dd{b}{t} & = & \tfrac{1}{2}\left(a-b\right)-\tfrac{1}{2}b^{2}+\dd{b}{x}+3\dd{^{2}b}{x^{2}}\,.
\end{eqnarray}
\end{subequations}
The lateral diffusion also included
in these \textsc{pde}s makes the derivation of macroscale boundary conditions challenging.

Various mathematical methodologies derive from \textsc{pde}s~\eqref{eq:heatExchanger} the macroscale model
\begin{equation}
\dd{C}{t}=\tfrac{1}{2}C^{3}-2C\dd{C}{x}+4\DD{C}{x}
+\Ord(C^{4}+\partial_{x}^{4}),\label{eq:slow}
\end{equation}
for the mean temperature~\(C(x,t):=[a(x,t)+b(x,t)]/2\)\,.
For example, centre manifold theory rigorously derives this effective macroscale model \citep{2013arXiv1310.1541R}, as does homogenization  \citep{PS08, MR2777986}.
This macroscale model combines an effective cubic reaction, with an effective nonlinear advection, and enhanced lateral diffusion.
The necessary analysis to derive the model~\eqref{eq:slow} is based around the equilibrium \(a=b=0\)\,, and applies to the slowly-varying in space solutions in the interior of the domain.
For example, in the interior it predicts the temperature fields are
\begin{equation}
\begin{bmatrix}a\\b\end{bmatrix}=
C\begin{bmatrix}1\\1\end{bmatrix}
+\left(\frac{1}{2}C^{2}-\dd{C}{x}\right)
\begin{bmatrix}1\\-1\end{bmatrix}
+\Ord(C^{3}+\partial_{x}^{3}).\label{eq:firstup}
\end{equation}
The challenge of this article is to provide sound boundary conditions for the macroscale model~\eqref{eq:slow}.

Prototypical microscale boundary conditions for microscale system~\eqref{eq:heatExchanger}
are taken to be the Dirichlet boundary conditions
\begin{equation}
a(0,t)=a_{0},\quad
b(0,t)=b_{0},\quad
a(L,t)=a_{L},\quad
\text{and}\quad b(L,t)=b_{L},\label{eq:microbc}
\end{equation}
where $a_{0}$, $a_{L}$, $b_{0}$ and~$b_{L}$ are potentially slowly varying functions of time. Section~\ref{sub:Projection-reveal-boundary} derives nonlinear macroscale boundary conditions~\eqref{eq:bc_constrans} for the macroscale
mean temperature model~\eqref{eq:slow}.
For example, the linearisation of the macroscale boundary condition~\eqref{eq:bc1} derived from the  Dirichlet boundary conditions~\eqref{eq:microbc} is the Robin condition
\begin{equation*}
C-\tfrac12\dd Cx\approx\tfrac14b_0+\tfrac34a_0\quad\text{at }x=0\,.
\end{equation*}
Importantly, this is not a Dirichlet boundary condition despite the microscale boundary conditions and the definition \(C:=(a+b)/2\) together suggesting boundary conditions are (incorrectly) \(C(0,t)=(a_0+b_0)/2\)\,.

\begin{figure}
 \centering\includegraphics{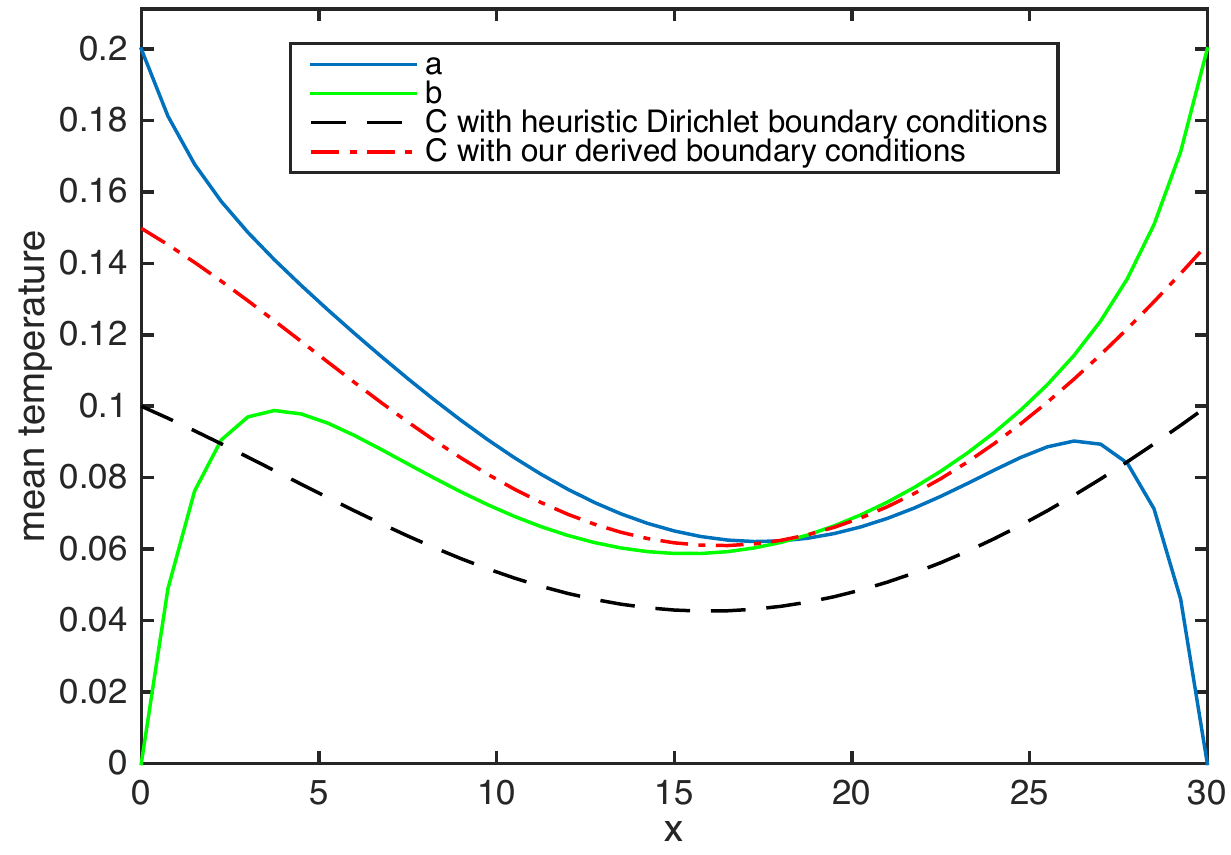}
 \caption{Example solutions of the heat exchanger~\eqref{eq:heatExchanger} in domain \(0\leq x\leq L=30\)
at time $t=21$. The two solid lines plot the the temperature of the two pipes, $a(x,t)$ and~$b(x,t)$. The dashed lines are solution of the macroscale
model~\eqref{eq:slow} at $t=21$\,: (black dashed)  with
heuristic Dirichlet boundary conditions;  and  (red dash-dots)  with
our systematically derived boundary conditions. \label{fig:boundarylayer} }
\end{figure}

Figure \ref{fig:boundarylayer} plots microscale and macroscale solutions
for the heat exchanger at a particular time.  The two solid lines plot the microscale solution~$a(x,t)$ and~$b(x,t)$ of microscale \textsc{pde}~\eqref{eq:heatExchanger}
 with microscale boundary conditions~\eqref{eq:microbc}. The black
dashed line plots the mean temperature model~\eqref{eq:slow} with classic Dirichlet boundary conditions $C(0,t)=\left(a_{0}+b_{0}\right)/2$ and $C(L,t)=\left(a_{L}+b_{L}\right)/2$ as would be commonly invoked \citep{MR2777986, multiscalethesis, Ray2012374}.
The macroscale model~\eqref{eq:slow} performs poorly with these heuristic Dirichlet boundary
conditions, especially in the interior of the domain (here \(5\leq x\leq 25\)). But the interior is where the macroscale model~\eqref{eq:slow} should be valid. The macroscale model~\eqref{eq:slow} represents
the interior dynamics but cannot resolve the details of boundary
layers \citep{Roberts92c}. With our derived boundary conditions,
the macroscale solution (red line in Figure~\ref{fig:boundarylayer})
fits the microscale solution (solid lines) in the interior:
the microscale fields \(a(x,t)\) and~\(b(x,t)\) being given by equation~\eqref{eq:firstup}.
Our systematic derivation of boundary conditions is needed for macroscale models to correctly predict the interior dynamics.

The key to our approach is to explore the effect of boundary layers by treating space as a time-like variable \citep[e.g.]{Chen2014}.
However, the heat exchanger problem~\eqref{eq:heatExchanger} is challenging because of the nonlinearity.
Here, a normal form coordinate transformation separates the spatial evolution in the boundary layers into a slow manifold, stable manifold and unstable manifolds.
This separation empowers a transformation of the given physical boundary conditions~\eqref{eq:microbc} into boundary conditions~\eqref{eq:bcgeneralformula} for the macroscale interior model~\eqref{eq:slow}.

\section{A normal form of the spatial evolution}
\label{sec:nfse}

The macroscale model~\eqref{eq:slow} is slow so the dominant terms in the boundary layers are due to the derivatives of spatial structure.
Thus, to derive macroscale boundary conditions for slow evolution~\eqref{eq:slow}
we treat the time derivative~$\partial/\partial t$ as a negligible operator~\citep{Roberts92c}.
To put heat exchanger system~\eqref{eq:heatExchanger} into the form of a dynamical system in time-like variable~\(x\) we define $a':=\dd{a}{x}$ and $b':=\dd{b}{x}$.
Then rearranging system~\eqref{eq:heatExchanger} in dynamical system form, with \(\partial_t=0\) for quasi-steady solution, gives
\begin{equation}
\dd{}{x}\begin{bmatrix}a\\
b\\
a'\\
b'
\end{bmatrix}=\begin{bmatrix}0 & 0 & 1 & 0\\
0 & 0 & 0 & 1\\
\frac{1}{6} & -\frac{1}{6} & \tfrac{1}{3} & 0\\
-\frac{1}{6} & \frac{1}{6} & 0 & -\tfrac{1}{3}
\end{bmatrix}\begin{bmatrix}a\\
b\\
a'\\
b'
\end{bmatrix}+\begin{bmatrix}0\\
0\\
-\frac{1}{2}a^{2}\\
\frac{1}{2}b^{2}
\end{bmatrix}.\label{eq:dcdiffusion}
\end{equation}
We analyse the (spatial) dynamics of this system with `initial condition' at \(x=0\) of the given microscale Dirichlet boundary conditions~\eqref{eq:microbc}.

Start by basing the analysis of~\eqref{eq:dcdiffusion} around the equilibrium at the origin, \(a=b=a'=b'=0\)\,.
The eigenvalues of the system linearised about the origin are $0$~(twice) and~\(\pm\sqrt{2}\)\,.
The eigenvalues of zero corresponds to an eigenvector of $\left(1,1,0,0\right)$
and a generalised eigenvector of $\left(-1,1,1,1\right)$. Hence the spatial \textsc{ode} system~\eqref{eq:dcdiffusion} contains two centre (slow) modes, one stable mode and one unstable mode.

\citet{WEB, MEDiCS} provides a web service
to construct by computer algebra a coordinate transform which separates stable, unstable and centre manifolds. However, the web service
does not directly apply to systems whose linearisation has a generalised eigenvector.
To circumvent the generalised eigenvector, we choose to
embed the \textsc{ode} system~\eqref{eq:dcdiffusion} as the \(\epsilon=1\) member of the one parameter family of systems
\begin{equation}
\dd{}{x}\begin{bmatrix}a\\
b\\
a'\\
b'
\end{bmatrix}=\begin{bmatrix}0 & 0 & 1 & -1\\
0 & 0 & -1 & 1\\
\frac{1}{6} & -\frac{1}{6} & -\frac{1}{6} & \frac{1}{2}\\
-\frac{1}{6} & \frac{1}{6} & -\frac{1}{2} & \frac{1}{6}
\end{bmatrix}\begin{bmatrix}a\\
b\\
a'\\
b'
\end{bmatrix}+\begin{bmatrix}\epsilon b'\\
\epsilon a'\\
-\tfrac{1}{2} a^{2}+\tfrac{1}{2}\epsilon a'-\tfrac{1}{2}\epsilon b'\\
\tfrac{1}{2} b^{2}+\tfrac{1}{2}\epsilon a'-\tfrac{1}{2}\epsilon b'
\end{bmatrix},\label{eq:dcdiffusion-1}
\end{equation}
where the last vector on the \textsc{rhs} is treated as a perturbative term: the parameter~$\epsilon$ counts the order of artificial linear perturbation.
The linear operator in system~\eqref{eq:dcdiffusion-1} now has no generalised eigenvector:
its eigenvalues are \(0\)~(twice) and~\(\pm\frac23\),
with corresponding eigenvectors $\left(1,1,0,0\right)$, $\left(-1,1,1,1\right)$, $\left(-\tfrac{3}{2},\tfrac{3}{2},0,1\right)$, $\left(-\tfrac{3}{2},\tfrac{3}{2},1,0\right)$.
The web service \citep{WEB} then finds a normal form coordinate transform as a multivariate power series in variables~\(s_j\) and parameter~$\epsilon$.
Substituting $\epsilon=1$ into the results reveals the centre manifold, stable manifold and unstable manifold for the spatial \textsc{ode}~\eqref{eq:dcdiffusion}.

The three manifolds can be parametrised as we choose.
We choose the definition of the two parameters
for the slow manifold to be the mean temperature $s_1=C$ and its spatial derivative $s_2=\dd{C}{x}$\,:
\begin{align} &
s_{1} := \tfrac{1}{2}\left(a+b\right)=C, &&
s_{3} := \tfrac{1}{8}\left(3a-3b-3a'+9b'\right),\nonumber \\ &
s_{2} := \tfrac{1}{2}\left(a'+b'\right)=\dd{C}{x}, &&
s_{4} := \tfrac{1}{8}\left(3a-3b+9a'-3b'\right),\label{eq:ctrans}
\end{align}
where $s_{3}$~parametrise the stable manifold, and $s_{4}$~parametrise the unstable manifold. Then the web service \citep{WEB} derives the coordinate transform~\eqref{eq:spatial_model} giving~$a$, \(b\), \(a'\)  and~$b'$ as a power series of~$s_{1}$, \(s_{2}\), \(s_{3}\), $s_{4}$, and~\(\epsilon\):
\begin{subequations}\label{eq:spatial_model}%
\begin{eqnarray}
a & \approx & s_{1}-s_{2}+0.25s_{3}+1.5s_{1}^{2}+6s_{2}^{2}-1.1s_{1}s_{3}-3.4s_{2}s_{3}-0.035s_{3}^{2}\nonumber \\
 &  & {}+0.75s_{4}+0.74s_{2}s_{4}+0.56s_{3}s_{4}-0.25s_{4}^{2}
 \,,\label{eq:spatial_modela} \\
b & \approx & s_{1}+s_{2}-0.75s_{3}-1.5s_{1}^{2}-6s_{2}^{2}-0.74s_{2}s_{3}+0.25s_{3}^{2}\nonumber \\
 &  & {}-0.25s_{4}-1.1s_{1}s_{4}+3.4s_{2}s_{4}-0.56s_{3}s_{4}-0.035s_{4}^{2}
 \,,\quad\label{eq:spatial_modelb} \\
a' & \approx & s_{2}+1.5s_{1}s_{2}-0.17s_{3}+0.56s_{1}s_{3}+0.91s_{2}s_{3}+0.47s_{3}^{2}\nonumber \\
 &  & {}+0.5s_{4}-0.56s_{1}s_{4}+1.2s_{2}s_{4}-0.33s_{4}^{2}
 \,,\label{eq:spatial_modelc} \\
b' & \approx & s_{2}-1.5s_{1}s_{2}+0.5s_{3}+0.56s_{1}s_{3}+1.2s_{2}s_{3}-0.33s_{3}^{2}\nonumber \\
 &  & {}-0.17s_{4}-0.56s_{1}s_{4}+0.91s_{2}s_{4}-0.047s_{4}^{2}
 \,,\label{eq:spatial_modeld}
\end{eqnarray}
\end{subequations}
For simplicity we only record these and later expressions correct to quadratic terms in~\(s_j\), that is, with cubic errors in the multinomial, and for simplicity we record coefficients to two significant figures, and here we evaluate the power series at \(\epsilon=1\) to recover a coordinate transform applicable to the original spatial system~\eqref{eq:dcdiffusion}.
The corresponding evolution of the spatial system~\eqref{eq:dcdiffusion} in these new variables~\(s_j\) is also provided by the web service which determines
\begin{subequations}\label{eq:spatial_evo}%
\begin{eqnarray}
\dd{s_{1}}{x} & \approx & s_{2}
\,,\label{eq:spatial_evoa} \\
\dd{s_{2}}{x} & \approx & 1.5s_{1}s_{2}
\,,\label{eq:spatial_evob} \\
\dd{s_{3}}{x} & \approx & -0.67s_{3}-0.75s_{3}s_{1}-0.94s_{3}s_{2}
\,,\label{eq:spatial_evoc} \\
\dd{s_{4}}{x} & \approx & +0.67s_{4}-0.75s_{4}s_{1}+0.94s_{4}s_{2}
\,.\label{eq:spatial_evod}
\end{eqnarray}
\end{subequations}
The normal form of the transformed system~\eqref{eq:spatial_evo} has useful properties.
Since $\dd{s_{3}}{x}=g_3(s_{1},s_{2})s_{3}$ and $\dd{s_{4}}{x}=g_4(s_{1},s_{2})s_{4}$ for some functions~$g_j$ indicates that three invariant manifolds of the system~\eqref{eq:spatial_evo} are \(s_3=0\)\,, \(s_4=0\) and \(s_3=s_4=0\).
From the linearisation of~\eqref{eq:spatial_evo} these are the centre-unstable, centre-stable, and slow manifolds respectively.
Further, because $\dd{s_{1}}{x}$ and~$\dd{s_{2}}{x}$ are functions of only~$s_{1}$ and~$s_{2}$, the planes of \(s_1\) and~\(s_2\) constant are isochrons of the slow manifold \citep{ANZ:3976228} (sometimes called the leaves of the foliation, {fibres}, a {fibration},
{fibre map}s or {fibre bundle}s \cite[pp.300--2, e.g.]{Murdock03}).

One might query whether the transformation~\eqref{eq:spatial_model} and~\eqref{eq:spatial_evo} is valid given that it is obtained by a power series in artificial parameter~\(\epsilon\) that is then evaluated at \(\epsilon=1\)\,.
The coefficients appear to converge well to the given values, but as an independent check we also embedded the spatial \textsc{ode}~\eqref{eq:dcdiffusion} into the different family
of problems
\begin{equation}
\dd{}{x}\begin{bmatrix}a\\
b\\
a'\\
b'
\end{bmatrix}=\begin{bmatrix}0 & 0 & 1 & -1\\
0 & 0 & -1 & 1\\
\frac{1}{6} & -\frac{1}{6} & \frac{1}{6} & \frac{1}{6}\\
-\frac{1}{6} & \frac{1}{6} & -\frac{1}{6} & -\frac{1}{6}
\end{bmatrix}\begin{bmatrix}a\\
b\\
a'\\
b'
\end{bmatrix}+\begin{bmatrix}\epsilon b'\\
\epsilon a'\\
-\tfrac{1}{2} a^{2}+\tfrac{1}{6}\epsilon a'-\tfrac{1}{6}\epsilon b'\\
\tfrac{1}{2} b^{2}+\tfrac{1}{6}\epsilon a'-\tfrac{1}{6}\epsilon b'
\end{bmatrix}.\label{eq:dcdiffusion-1-1}
\end{equation}
Performing the same algebraic construction, but from this quite different base, we find
system~\eqref{eq:dcdiffusion-1-1} results in the same transform~\eqref{eq:spatial_model} and evolution~\eqref{eq:spatial_evo}.
This confirms the perturbative approach via embedding.

\section{Projection reveals boundary conditions\label{sub:Projection-reveal-boundary}}

This section focuses on the boundary layer near \(x=0\). As shown by the solid lines in Figure~\ref{fig:boundarylayer}, the microscale boundary conditions at $x=0$ force a boundary layer in the microscale model \eqref{eq:heatExchanger}. However, the macroscale model \eqref{eq:slow} does not resolve the boundary layer.
Forcing the macroscale model to pass through $(a_0+b_0)/2$  introduces an error in the interior of the domain, as shown by the dashed blue line in Figure~\ref{fig:boundarylayer}.
Here we derive an improved boundary condition at $x=0$ which reduces the interior error caused by poorly chosen macroscale boundary condition.

The boundary layer must lie in the centre-stable manifold \(s_4=0\) because if there was any component \(s_4\neq0\) then this would grow exponentially quickly in space and dominate the solution across the whole domain.
Algebraically we obtain the centre-stable manifold by substituting \(s_4=0\) into the coordinate transform~\eqref{eq:spatial_model}: the terms in~\eqref{eq:spatial_model} are arranged so that this simply means omitting the second line of each of the four pairs of lines.

\begin{figure}
\centering\includegraphics{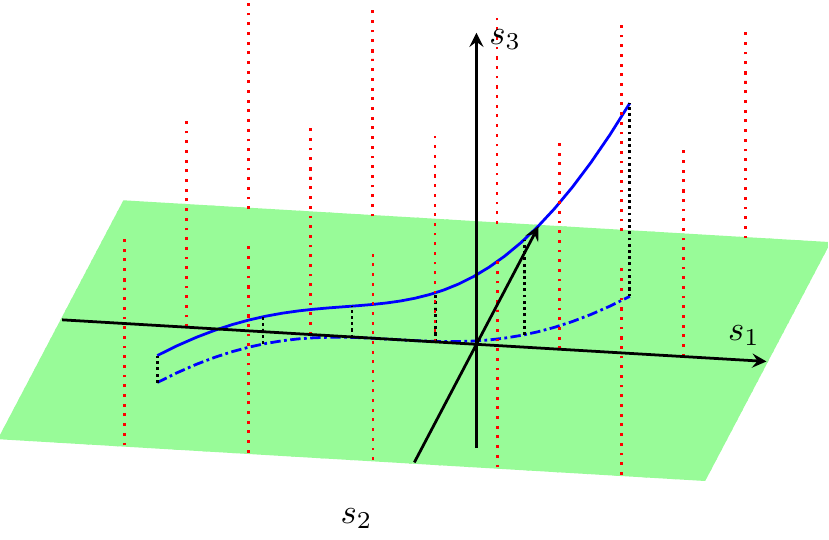}
\caption{schematic plot of centre-stable manifold near the boundary $x=0$. The green
plane is the centre manifold. The blue solid line is the set of values allowed by the microscale
boundary condition at $x=0$. The blue dotted line is the projection
of the microscale boundary values onto the slow manifold. \label{fig:cenre-stable} }
\end{figure}

Then, as plotted schematically in Figure~\ref{fig:cenre-stable}, the two Dirichlet boundary conditions~\eqref{eq:microbc} at $x=0$ form a one dimensional
curve (solid blue line)  of allowed values in the three-dimensional centre-stable manifold parametrised by~$s_{1}$, $s_{2}$ and~$s_{3}$.
Recall $a_{0}$
and~$b_{0}$ are the boundary values at $x=0$ from boundary conditions~\eqref{eq:microbc}.
The first two components
on the centre-stable manifold (\(s_4=0\)) of~\eqref{eq:spatial_model} reveal the microscale constraints on the boundary, upon defining  $s_{i}^{0}:=s_{i}\big|_{x=0}$ for $i=1,2,3$\,,
\begin{equation}
\begin{bmatrix}a_{0}\\
b_{0}
\end{bmatrix}\approx\begin{bmatrix}s_{1}^{0}-s_{2}^{0}+0.25s_{3}^{0}+1.5s_{1}^{0}{}^{2}+6{{{s_{2}^{0}}}}^{2}-1.1s_{1}^{0}s_{3}^{0}-3.4s_{2}^{0}s_{3}^{0}-0.035{{{s_{3}^{0}}}}^{2}\\
s_{1}^{0}+s_{2}^{0}-0.75s_{3}^{0}-1.5{{{s_{1}^{0}}}}^{2}-6{s_{2}^{0}}^{2}-0.74s_{2}^{0}s_{3}^{0}+0.25{{{s_{3}^{0}}}}^{2}
\end{bmatrix}.\label{eq:evo3}
\end{equation}
These equations implicitly determines the solid blue curve in Figure~\ref{fig:cenre-stable}.
To explicitly describe the curve, recall that this is a power series with cubic errors and so we just need to consistently revert the series to give, say, the boundary values~\(s_1^0\) and~\(s_3^0\) as a function of~\(s_2^0\), \(a_0\) and~\(b_0\).
Algebra determines
\begin{subequations}\label{eq:bc_constrans}%
\begin{eqnarray}
s_{1}^{0} & \approx & \left(0.25b_{0}-0.29b_{0}^{2}+0.75a_{0}-0.63a_{0}b_{0}+0.18a_{0}^{2}\right)
\nonumber\\&&{}
+s_{2}^{0}\left(0.5-2.8b_{0}+3.7a_{0}\right)+3{s_{2}^{0}}^{2},\label{eq:bc_constransa} \\
s_{3}^{0} & \approx & \left(-b_{0}-0.19b_{0}^{2}+a_{0}-2.3a_{0}b_{0}-0.56a_{0}^{2}\right)
\nonumber\\&&{}
+s_{2}^{0}\left(2-4.6b_{0}+3.8a_{0}\right)-5.2{s_{2}^{0}}^{2}.\label{eq:bc_constransb}
\end{eqnarray}
\end{subequations}

Since the slow dynamics in the interior of the domain must lie on the slow manifold \(s_3=0\), appropriate boundary conditions for the interior dynamics must come from projecting these allowed boundary values onto the slow manifold.
Because of the special normal form of the transformed system~\eqref{eq:spatial_evo},  the slow variables \(s_1\) and~\(s_2\) evolve independently of the fast variables~\(s_3\) and~\(s_4\), and the appropriate projection is the orthogonal projection along the isochrons \(s_1\) and~\(s_2\) constant onto the plane \(s_3=0\) ---shown by the red lines in Figure~\ref{fig:cenre-stable}.
Equation~\eqref{eq:bc_constransa} describes the projected curve in the \(s_1s_2\)-plane illustrated by the blue dashed line in Figure~\ref{fig:cenre-stable}. Recall from the amplitude definition \eqref{eq:ctrans} that $C$  and $s_1$ are the same. Hence substituting $s_1^0=C$ and $s_2^0=\frac{\partial C}{\partial x}$ into equation~\eqref{eq:bc_constransa} forms the boundary condition at $x=0$
\begin{eqnarray}&&
C-\left(0.5-2.8b_{0}+3.7a_{0}\right)\dd{C}{x}-3\left(\dd{C}{x}\right)^{2}
\nonumber\\&&
\approx\left(0.25b_{0}-0.29b_{0}^{2}+0.75a_{0}-0.63a_{0}b_{0}+0.18a_{0}^{2}\right).
\label{eq:bcgeneralformula}
\end{eqnarray}
This nonlinear Robin boundary condition produces the correct macroscale slowly varying interior domain solutions of the microscale model \textsc{pde}~\eqref{eq:heatExchanger}.

\section{A numerical example}

As an example, let the boundary values be $a_{0}=0.2f(t)$ and $b_{0}=0$ for \(f(t)=\tanh^2t\) varying smoothly but quickly from \(f(0)=0\) to $1$.  
 Macroscale boundary condition~\eqref{eq:bcgeneralformula}
gives the macroscale boundary condition at $x=0$ for mean temperature
model~\eqref{eq:slow}
\begin{equation}
C-\left[0.75f + 0.5\right]\dd{C}{x}-3\left(\dd{C}{x}\right)^{2}=0.15f + 0.007f^2.\label{eq:bc1}
\end{equation}

%

\paragraph{Macroscale boundary conditions on the right}

One method to  derive the macroscale boundary conditions
at $x=L$ is to appeal to symmetry.
Define a new spatial coordinate $\tilde x=L-x$ measuring distance from the boundary into the interior, and define new field variables \(\tilde a(\tilde x,t)=-b(x,t)\), \(\tilde b(\tilde x,t)=-a(x,t)\) and therefore \(\tilde C(\tilde x,t)=-C(x,t)\).
Then the \textsc{pde} system~\eqref{eq:heatExchanger} is symbolically identical in the tilde and plain variables.
But the boundary conditions~\eqref{eq:microbc} at the right-boundary \(x=L\) are transformed to Dirichlet boundary conditions at \(\tilde x=0\) of \(\tilde a(0,t)=-b_L\) and \(\tilde b(0,t)=-a_L\)\,.
Then the derivation of Sections~\ref{sec:nfse} and~\ref{sub:Projection-reveal-boundary} apply in the same way to the tilde problem.
After computing the macroscale boundary conditions in coordinate~$\tilde x$ we transform back to the original coordinate~$x$.

For example, assume $a_{L}=0$ and $b_{L}=0.2$.
The iteration scheme in Section~\ref{sub:Projection-reveal-boundary} computes macroscale boundary
condition on the boundary $x=L$ ($\tilde x=0$)
\begin{equation}
-\tilde C-\left[0.75f - 0.5\right]\dd{\tilde C}{\tilde x}+3\left(\dd{\tilde C}{\tilde x}\right)^{2}=0.15f -0.007f^2.\label{eq:paper2bc2-0}
\end{equation}
By the chain rule $\dd{\tilde x}{x}=-1$\,, and substitute \(\tilde C(\tilde x,t)=-C(x,t)\) into boundary condition \eqref{eq:paper2bc2-0}
\begin{equation}
C-\left[0.75f - 0.5\right]\dd{C}{x}+3\left(\dd{C}{x}\right)^{2}=0.15f -0.007f^2.\label{eq:bc2-1}
\end{equation}

\paragraph{Numerics verifies the macroscale boundary conditions derivation}

Figure \ref{fig:boundarylayer} plots a snapshot of the simulations
on microscale model~\eqref{eq:heatExchanger} and mean temperature
model~\eqref{eq:slow} for two cases: the Dirichlet boundary conditions $C_{0}=\left(a_{0}+b_{0}\right)/2$
and $C_{L}\left(a_{L}+b_{L}\right)/2$; and our systematic boundary conditions~\eqref{eq:bc1} and~\eqref{eq:bc2-1}.
Using finite differences we
convert the system of two \textsc{pde}s~\eqref{eq:heatExchanger} into a system of \textsc{ode}s.
Then Matlab's \texttt{ode15s}  applies a variable
order method to compute the solution of the system of \textsc{ode}s  \citep{Shampine:1999:SID:333769.333775}.

The numerical result is as expected.
The macroscale model with systematic boundary conditions~\eqref{eq:bcgeneralformula} model the interior domain microscale dynamics much better than that with heuristic Dirichlet boundary conditions.

\section{Conclusion\label{sec:con-1}}

We systematically derived macroscale boundary conditions from microscale
Dirichlet boundary conditions. This methodology can be extended to
microscale Neumann and Robin boundary conditions. For the microscale Dirichlet boundary
conditions, we evaluated the first two components of the centre-stable
manifold~\eqref{eq:spatial_modela}--\eqref{eq:spatial_modelb} at $x=0$ to reveal the microscale boundary
constraints~\eqref{eq:evo3}.
If the microscale boundary conditions
were Neumann, we would use the last two components,~\eqref{eq:spatial_modelc}--\eqref{eq:spatial_modeld}. If the microscale
boundary conditions were Robin, we would use linear combinations of
the transform~\eqref{eq:spatial_model}.
The methodology also applies to more general multiscale modelling of \textsc{pde}s \citep{Roberts92c}.

\paragraph{Acknowledgements}
CC thanks Dr.~Tony Miller for his advice and useful discussion, and
 \textsc{csiro} for their support in funding to participate in conferences and workshops.

\bibliographystyle{agsm}
\bibliography{research}

\end{document}